\documentclass[12pt]{amsart}
\usepackage{amssymb,amsmath}

\numberwithin{equation}{section}

\newtheorem{Thm}{Theorem}[section]
\newtheorem{Prop}{Proposition}[section]

\newtheorem{Lem}{Lemma}[section]

\def \eps{{\epsilon}}

\def \EE{{\mathcal E}}

\def \KK{{\mathcal K}}

\def \OO{{\mathcal O}}

\def \QQ{{\mathcal Q}}

\def \gHat{{\widehat g}}

\def \gTld{{\tilde g}}

\def \One{{\mathbf{1}}}

\def \ini{\mathrm{in}}

\def \inttwo{\int \!\!\! \int }
\def \intthree{\int \!\!\! \int \!\!\! \int }

\def \dee{{\mathrm d}}
\def \dmu{{\dee \mu}}
\def \dsig{{\dee \sigma_{\! x}}}
\def \domega{{\dee \omega}}

\def \ds{{\dee s}}
\def \dt{{\dee t}}
\def \dx{{\dee x}}
\def \dv{{\dee v}}

\def \del{{\partial}}
\def \GRAD{\nabla_{\!\!x}}

\def \DIV{\nabla_{\!\!x} \! \cdot }

\def \DOT{{\,\cdot\,}}

\def \BAR{\overline}

\def \<{\langle}
\def \>{\rangle}

\def \lANGLE{\big\langle \! \! \big\langle}
\def \rANGLE{\big\rangle \! \! \big\rangle}
\def \LANGLE{\bigg\langle \! \! \! \! \! \; \bigg\langle}
\def \RANGLE{\bigg\rangle \! \! \! \! \! \; \bigg\rangle}
\def \p{\partial}
\def \grad{\nabla_{\!x}}
\def \pO{\partial\Omega}

\def \gpg{\gamma_+g_\epsilon}
\def \gpgb{\langle\gamma_+g_\epsilon\rangle_{\partial\Omega}}
\def \gamp{\gamma_+}
\def \gamn{\gamma_-}

\def\R{{\mathbb R}}      

\def\D{{\mathrm{D}}}

\def\RD{{{\mathbb R}^{\D}}}
\def\SD{{{\mathbb S}^{\D-1}}}

\def\RRDR{{\R \times \R^\D \! \times \R}}

\def\SDRD{{\SD \! \times \RD}}
\def\SDRDRD{{\SD \! \times \RD \! \times \RD}}

\def \wL{\hbox{w-}L}

\def \init{{\mathrm{in}}}

\def \dep{\delta_\eps}
\def \gps{\tilde{g}_\eps}
\def \Ge{G_\eps}

\def \bBar{\BAR{b}}

\def\nn{\mathrm{n}}
\def\Rn{\mathrm{R}}

\def\nRm{\mathrm{n}}
\def\RRm{\mathrm{R}}

\begin{document}

\title[Acoustic Limit for the Boltzmann Equation]
      {Remarks on the Acoustic Limit for \\ the Boltzmann Equation}

\author[Ning Jiang]{Ning Jiang}
\address[N.J.]%
        {Courant Institute of Mathematical Sciences,
         251 Mercer Street, New York, NY 10012}
\email{njiang@cims.nyu.edu}
\author[C. David Levermore]{C. David Levermore}
\address[C.D.L.]%
        {Department of Mathematics {\em and}
         Institute for Physical Science and Technology,
         University of Maryland, College Park, MD 20742}
\email{lvrmr@math.umd.edu}
\author[Nader Masmoudi]{Nader Masmoudi}
\address[N.M.]%
        {Courant Institute of Mathematical Sciences,
         251 Mercer Street, New York, NY 10012}
\email{masmoudi@cims.nyu.edu}

\date{\today}

\begin{abstract}

   We improve in three ways the results of \cite{GL} that establish 
the acoustic limit for DiPerna-Lions solutions of Boltzmann equation.  
First, we enlarge the class of collision kernels treated to that 
found in \cite{LM}, thereby treating all classical collision kernels 
to which the DiPerna-Lions theory applies.  Second, we improve the 
scaling of the kinetic density fluctuations with Knudsen number from
$O(\epsilon^m)$ for some $m>\frac12$ to $O(\epsilon^\frac12)$.  Third, 
we extend the results from periodic domains to bounded domains with 
a Maxwell reflection boundary condition, deriving the impermeable 
boundary condition for the acoustic system.

\end{abstract}

\maketitle


\section{Introduction}
In this note we establish the acoustic limit starting from
DiPerna-Lions renormalized solutions of the Boltzmann equation
considered over a smooth bounded spatial domain $\Omega\subset\RD$.
The acoustic system is the linearization about the homogeneous state
of the compressible Euler system.  After a suitable choice of units
and Galilean frame, it governs the fluctuations in mass density
$\rho(x,t)$, bulk velocity $u(x,t)$, and temperature $\theta(x,t)$
over $\Omega\times\R_+$ by the initial-value problem
\begin{equation}
  \label{Acoustic-IVP}
\begin{aligned}
  \del_t \rho + \DIV u & = 0 \,, & \qquad
  \rho(x,0) & = \rho^\init(x) \,,
\\
  \del_t u + \GRAD (\rho + \theta) & = 0 \,, & \qquad
  u(x,0) & = u^\init(x) \,,
\\
  \tfrac{\D}{2} \del_t \theta + \DIV u & = 0 \,, & \qquad
  \theta(x,0) & = \theta^\init(x) \,,
\end{aligned}
\end{equation}
subject to the impermeable boundary condition
\begin{equation}
  \label{Acoustic-BC}
  u \DOT \nRm = 0 \,, \quad \mbox{on $\del\Omega$} \,,
\end{equation}
where $\nRm(x)$ is the unit outward normal at $x\in\del\Omega$.  
This is one of the simplest fluid dynamical systems imaginable, 
being essentially the wave equation.

   The acoustic system (\ref{Acoustic-IVP}, \ref{Acoustic-BC}) can
be formally derived from the Boltzmann equation for kinetic densities
$F(v,x,t)$ over $\RD\times\Omega\times\R_+$ that are close to the
global Maxwellian
\begin{equation}
  \label{abs-Max}
  M(v) = \frac{1}{(2\pi)^{\frac{\D}2}}
         \exp\!\big( - \tfrac12 |v|^2 \big) \,.
\end{equation}
We consider families of kinetic densities in the form
$F_\eps(v,x,t)=M(v)G_\eps(v,x,t)$ where the relative kinetic
densities $G_\eps(v,x,t)$ over $\RD\times\Omega\times\R_+$ are
governed by the rescaled Boltzmann initial-value problem
\begin{equation}
  \label{Boltzmann-IVP}
  \del_t G_\eps + v \DOT \GRAD G_\eps
  = \frac{1}{\eps} \QQ(G_\eps,G_\eps) \,, \qquad
  G_\eps(v,x,0) = G_\eps^\init(v,x) \,.
\end{equation}
Here the Knudsen number $\eps>0$ is the ratio of the mean free path
to a macroscopic length scale and the collision operator
$\QQ(G_\eps,G_\eps)$ is given by
\begin{equation}
  \label{Coll-Op}
  \QQ(G_\eps,G_\eps)
  = \inttwo_\SDRD
           \big( G'_{\eps 1} G'_\eps - G_{\eps 1} G_\eps \big) \,
           b(\omega,v_1-v) \, \domega \, M_1 \dv_1 \,,
\end{equation}
where the collision kernel $b(\omega,v_1-v)$ is positive almost
everywhere while $G_{\eps 1}$, $G'_\eps$, and $G'_{\eps 1}$ denote
$G_\eps(\,\cdot\,,x,t)$ evaluated at $v_1$,
$v'=v+\omega\omega\DOT(v_1-v)$, and $v'_1=v-\omega\omega\DOT(v_1-v)$
respectively.

   We impose a Maxwell reflection boundary condition on $\del\Omega$
of the form
\begin{equation}
  \label{Boltzmann-BC}
  \One_{\Sigma_+} G_\eps \circ \RRm
  = (1 - \alpha) \, \One_{\Sigma_+} \, G_\eps
    + \alpha \, \One_{\Sigma_+} \, \sqrt{2 \pi} \,
      \big\< \One_{\Sigma_+} \, v \DOT n \, G_\eps \big\> \,.
\end{equation}
Here $\alpha\in[0,1]$ is the Maxwell accommodation coefficient,
$(G_\eps\circ\RRm)(v,x,t)=G_\eps(\RRm(x)v,x,t)$ where
$\RRm(x)=I-2\nRm(x)\nRm(x)^T$ is the specular reflection matrix at a
point $x\in\del\Omega$, $\One_{\Sigma_+}$ is the indicator function
of the so-called outgoing boundary set
\begin{equation}
  \label{Outgoing-Boundary}
  \Sigma_+ = \left\{ (v,x) \in \RD \times \del\Omega \,:\,
                     v \DOT \nRm(x) > 0 \right\} \,,
\end{equation}
and $\<\,\cdot\,\>$ denotes the average
\begin{equation}
  \label{Maxwell-Avg}
  \< \xi \> = \int_\RD \xi(v) \, M(v) \, \dv \,.
\end{equation}
Because $\sqrt{2\pi}\,\big\<\One_{\Sigma_+}\,v\DOT\nRm\big\>=1$, it 
is easy to see from \eqref{Boltzmann-BC} that on $\del\Omega$ one has
\begin{equation}
  \label{Boundary-Flux}
\begin{aligned}
  \< v \DOT \nRm \, G_\eps \>
  & = \big\< \One_{\Sigma_+} \, v \DOT \nRm \,
             \big( G_\eps - G_\eps \circ \Rn \big) \big\>
\\
  & = \alpha \,
      \Big\< \One_{\Sigma_+} \, v \DOT \nRm \,
             \Big( G_\eps
                   - \sqrt{2 \pi} \,
                     \big\< \One_{\Sigma_+} \, v \DOT \nRm \,
                            G_\eps \big\> \Big) \Big\> = 0 \,.
\end{aligned}
\end{equation}

   Fluid regimes are those in which the Knudsen number $\eps$ is
small.  The acoustic system (\ref{Acoustic-IVP}, \ref{Acoustic-BC})
can be derived from (\ref{Boltzmann-IVP}, \ref{Boltzmann-BC}) for
families of solutions $G_\eps(v,x,t)$ that are scaled so that
\begin{equation}
  \label{G-Fluct}
  G_\eps = 1 + \delta_\eps g_\eps \,, \qquad
  G_\eps^\init = 1 + \delta_\eps g_\eps^\init \,,
\end{equation}
where
\begin{equation}
  \label{delta-optimal}
  \delta_\eps \to 0 \qquad \hbox{as} \qquad \eps \to 0 \,,
\end{equation}
and the fluctuations $g_\eps$ and $g_\eps^\init$ converge in the
sense of distributions to $g\in L^\infty(\dt;L^2(M\dv\,\dx))$
and $g^\init\in L^2(M\dv\,\dx)$ respectively as $\eps\to0$.
One finds that $g$ has the infinitesimal Maxwellian form
\begin{equation}
  \label{Inf-Max}
  g = \rho + v \DOT u
           + \big( \tfrac12 |v|^2 - \tfrac{\D}2 \big) \theta \,,
\end{equation}
where $(\rho,u,\theta)\in L^\infty(\dt;L^2(\dx;\RRDR))$ solve
(\ref{Acoustic-IVP}, \ref{Acoustic-BC}) with initial data given by
\begin{equation}
  \label{Init-Data}
  \rho^\init = \< g^\init \> \,, \qquad
  u^\init = \< v \, g^\init \> \,, \qquad
  \theta^\init = \big\< \big( \tfrac1\D |v|^2 - 1 \big) \,
                        g^\init \big\> \,.
\end{equation}
The formal derivation leading to \eqref{Acoustic-IVP} closely
follows that in \cite{GL}, so its details will not be given here.
The boundary condition \eqref{Acoustic-BC} is obtained by noticing
that \eqref{Boundary-Flux} implies $\<v\DOT\nRm\,g_\eps\>=0$, then
passing to the limit in this to get $\<v\DOT\nRm\,g\>=0$, and
finally using \eqref{Inf-Max} to obtain \eqref{Acoustic-BC}.

   The program initiated in \cite{BGL1, BGL2, BGL3} seeks to justify
fluid dynamical limits for Boltzmann equations in the setting of
DiPerna-Lions renormalized solutions \cite{DiPL}, which are the
only temporally global, large data solutions available.  The main
obstruction to carrying out this program is that DiPerna-Lions
solutions are not known to satisfy many properties that one formally
expects for solutions of the Boltzmann equation.  For example, they
are not known to satisfy the formally expected local conservations
laws of momentum and energy.  Moreover, their regularity is poor.
The justification of fluid dynamical limits in this setting is
therefore not easy.

   The acoustic limit was first established in this kind of setting
in \cite{BGL} over a periodic domain.  There idea introduced there
was to pass to the limit in approximate local conservations laws which
are satified by DiPerna-Lions solutions.  One then shows that the
so-called conservation defects vanish as the Knudsen number $\eps$
vanishes, thereby establishing the local conservation laws in the
limit.  This was done in \cite{BGL} using only relative entropy
estimates, which restricted the result to collision kernels that are
bounded and to fluctuations scaled so that
\begin{equation}
  \label{BGL-scaling}
  \delta_\eps \to 0 \quad \hbox{and} \quad
  \frac{\delta_\eps}{\eps} |\log(\delta_\eps)| \to 0
  \quad \hbox{as} \quad \eps \to 0 \,,
\end{equation}
which is far from the formally expected optimal scaling
\eqref{delta-optimal}.

   In \cite{GL} the local conservation defects were removed using new
dissipation rate estimates.  This allowed the treatment of collision
kernels that for some $C_b<\infty$ and $\beta\in[0,1)$ satisfied
\begin{equation}
  \label{GL-kernel}
  \int_\SD b(\omega,v_1-v) \, \domega
  \leq C_b \big( 1 + |v_1 - v|^2 \big)^\beta \,,
\end{equation}
and of fluctuations scaled so that
\begin{equation}
  \label{GL-scaling}
  \delta_\eps \to 0 \quad \hbox{and} \quad
  \frac{\delta_\eps}{\eps^{1/2}} |\log(\delta_\eps)|^{\beta/2} \to 0
  \quad \hbox{as} \quad \eps \to 0 \,.
\end{equation}
The above class of collision kernels includes all classical kernels
that are derived from Maxwell or hard potentials and that satisfy a
weak small deflection cutoff.  The scaling given by (\ref{GL-scaling})
is much less restrictive than that given by (\ref{GL-scaling}), but is
far from the formally expected optimal scaling \eqref{delta-optimal}.
Finally, only periodic domains are treated in \cite{GL}.

   Here we improve the result of \cite{GL} in three ways.  First, we
apply estimates from \cite{LM} to treat a broader class of collision 
kernels that includes those derived from soft potentials.  Second, we 
improve the scaling of the fluctuations to $\delta_\eps=O(\eps^{1/2})$.  
Finally, we treat domains with a boundary and use new estimates to 
derive the boundary condition \eqref{Acoustic-BC} in the limit.

   We use the $L^1$ velocity averaging theory of Golse and
Saint-Raymond \cite{GS-1} through the nonlinear compactness estimate
of \cite{LM} to improve the scaling of the fluctuations to
$\delta_\eps=O(\eps^{1/2})$.  Without it we would only be
able to improve the scaling to $\delta_\eps=o(\eps^{1/2})$.
This is the first time the $L^1$ averaging theory has played any role
in an acoustic limit theorem, albeit for a modest improvement in the 
scaling of our result.  We remark that averaging theory plays no role 
in establishing the Stokes limit with its formally expected optimal 
scaling of $\delta_\eps=o(\eps)$ \cite{LM}.

   We treat domains with boundary in the setting of Mischler
\cite{Misch}, who extended the DiPerna-Lions theory to bounded
domains with a Maxwell reflection boundary condition.  He showed
that these boundary conditions are satisfied in a {\em renormalized}
sense.  This means we cannot deduce that $\<v\DOT\nRm\,g_\eps\>\to0$
as $\eps\to0$ to derive the boundary condition \eqref{Acoustic-BC},
as we did in our formal argument.  Masmoudi and Saint-Raymond 
\cite{MS} developed estimates to obtain boundary conditions in the 
Stokes limit.  However neither these estimates nor their recent 
extension to the Navier-Stokes limit \cite{JM} can handle the 
acoustic limit.  Rather, we develop new boundary {\em a priori} 
estimates to obtain a weak form of the boundary condition 
\eqref{Acoustic-BC} in this limit.  In doing so, we treat a 
broader class of collision kernels than was done in \cite{MS}.

   Finally, we remark that fully establishing the acoustic limit
with its formally expected optimal scaling of the fluctuation size
\eqref{delta-optimal} is still open.  This gap must be bridged
before one can hope to fully establish the compressible Euler limit
starting from DiPerna-Lions solutions to the Boltzmann equation. In
contrast, optimal scaling can be obtained within the framework of
classical solutions by using the nonlinear energy method developed
by Guo. This has been done recently by the first author of this
paper with Guo and Jang \cite{GJJ-1, GJJ-2}.

   Our paper is organized as follows.  Section 2 gives its framework.
Section 3 states and proves our main result modulo two steps.
Section 4 removes the conservation defects.  Section 5 establishes
the limit boundary mass-flux term.


\section{Framework}

   For the most part we will use the notation of \cite{LM}.  Here we 
present only what is needed to state our theorem.  For more complete 
introductions to the Boltzmann equation, see \cite{GL,LM}.

   Let $\Omega\subset\RD$ be a bounded domain with smooth boundary
$\del\Omega$.  Let $\nRm(x)$ denote the outward unit normal vector at 
$x\in\del\Omega$ and $\dee\sigma_x$ denote the Lebesgue measure on 
$\del\Omega$.  The phase space domain associated with $\Omega$ is 
$\OO=\RD\times\Omega$, which has boundary 
$\del\OO=\RD\times\del\Omega$.  Let $\Sigma_+$ and $\Sigma_{-}$
denote the outgoing and incoming subsets of $\del\OO$ defined by
\begin{equation}
  \nonumber
  \Sigma_{\pm} = \left\{ (v,x) \in \del\OO \,:\,
                         \pm v \DOT \nRm(x) > 0 \right\} \,.
\end{equation}
The global Maxwellian $M(v)$ given by \eqref{abs-Max} corresponds
to the spatially homogeneous fluid state with density and temperature
equal to 1 and bulk velocity equal to 0.  The boundary condition
\eqref{Boltzmann-BC} corresponds to a wall temperature of 1, so that
$M(v)$ is the unique equilibrium of the fluid.  Associated with the 
initial data $G_\eps^\init$ we have the normalization
\begin{equation}
  \label{norm-init}
  \int_\Omega \< G_\eps^\init \> \, \dx = 1 \,.
\end{equation}


\subsection{Assumptions on the Collision Kernel}
\label{AssumpKernel}

   The kernel $b(\omega,v_1-v)$ associated with the collision operator
\eqref{Coll-Op} is positive almost everywhere.  The Galilean
invariance of the collisional physics implies that $b$ has the
classical form
\begin{equation}
  \label{KernelForm}
  b(\omega,v_1 - v)
  = |v_1 - v| \, \Sigma(|\omega \DOT n|,|v_1 - v|) \,,
\end{equation}
where $n=(v_1-v)/|v_1-v|$ and $\Sigma$ is the specific differential
cross-section.  We make five additional technical assumptions 
regarding $b$ that are adopted from \cite{LM}.

\smallskip

   Our {\em first technical assumption} is that the collision kernel
$b$ satisfies the requirements of the DiPerna-Lions theory.  That
theory requires that $b$ be locally integrable with respect to
$\domega\,M_1\dv_1\,M\dv$, and that it moreover satisfies
\begin{equation}
  \label{bLimit}
  \lim_{|v|\to\infty}
      \frac{1}{1 + |v|^2} \int_K \! \bBar(v_1 - v) \, \dv_1
  = 0 \quad
  \hbox{for every compact $K\subset\RD$} \,,
\end{equation}
where $\bBar$ is defined by
\begin{equation}
  \label{bBarDef}
  \bBar(v_1 - v) \equiv \int_\SD \! b(\omega,v_1 - v) \, \domega \,.
\end{equation}
Galilean symmetry (\ref{KernelForm}) implies that $\bBar$ is a
function of $|v_1 - v|$ only.

\smallskip

   Our {\em second technical assumption} regarding $b$ is that the
attenuation coefficient $a$, which is defined by
\begin{equation}
  \label{aDef}
  a(v) \equiv \int_\RD \! \bBar(v_1 - v) \, M_1 \dv_1
       = \inttwo_\SDRD \!\! b(\omega,v_1 - v) \, \domega \, M_1 \dv_1 \,,
\end{equation}
is bounded below as
\begin{equation}
  \label{aLowerBnd}
  C_a \big( 1 + |v|^2 \big)^{\beta_a} \leq a(v) \quad
  \hbox{for some constants $C_a>0$ and $\beta_a\in\R$} \,.
\end{equation}
Galilean symmetry (\ref{KernelForm}) implies that $a$ is a function 
of $|v|$ only.

\smallskip

   Our {\em third technical assumption} regarding $b$ is that there
exists $s\in(1,\infty]$ and $C_b\in(0,\infty)$ such that
\begin{equation}
  \label{bBnd}
  \left( \int_\RD \bigg| \frac{\bBar(v_1 - v)}
                              {a(v_1) \, a(v)} \bigg|^s
                  a(v_1) \, M_1 \dv_1 \right)^{\frac1s}
  \leq C_b \,.
\end{equation}
Because this bound is uniform in $v$, we may take $C_b$ to be the
supremum over $v$ of the left-hand side of (\ref{bBnd}).

\smallskip

   Our {\em fourth technical assumption} regarding $b$ is that the
operator
\begin{equation}
  \label{KK^+Cpt}
  \KK^+ : L^2(aM\dv) \to L^2(aM\dv) \quad
  \hbox{is compact} \,,
\end{equation}
where
\begin{equation}
  \nonumber
  \KK^+ \tilde{g}
  = \frac1{2a} \inttwo_\SDRD \!
                   \big( \tilde{g}' + \tilde{g}'_1 \big) \,
                   b(\omega,v_1 - v) \, \domega \, M_1 \dv_1 \,.
\end{equation}
We remark that $\KK^+:L^2(aM\dv)\to L^2(aM\dv)$ is always bounded
\cite{LM,LS} with $\|\KK^+\|\leq1$.

\smallskip

   Our {\em fifth technical assumption} regarding $b$ is that for
every $\delta>0$ there exists $C_\delta$ such that $\bBar$ satisfies
\begin{equation}
  \label{Sat-Bnd}
  \frac{\bBar(v_1 - v)}
       {1 + \delta \, \dfrac{\bBar(v_1 - v)}{1 + |v_1 - v|^2}}
  \leq C_\delta \big( 1 + a(v_1) \big)
                \big( 1 + a(v) \big) \quad
  \hbox{for every $v_1, v\in\RD$} \,.
\end{equation}

\smallskip

   The above assumptions are satisfied by all the classical collision
kernels with a weak small deflection cutoff that derive from a
repulsive intermolecular potential of the form $c/r^k$ with
$k>2\tfrac{\D-1}{\D+1}$.  This includes all the classical collision
kernels to which the DiPerna-Lions theory applies \cite{LM,LS}.  
Kernels that satisfy \eqref{GL-kernel} clearly satisfy \eqref{bLimit}.  
If they moreover satisfy \eqref{aLowerBnd} with $\beta_a=\beta$ then 
they also satisfy \eqref{bBnd} and \eqref{Sat-Bnd}.  

   Because the kernel $b$ satisfies \eqref{bLimit}, it can be
normalized so that
\begin{equation}
  \nonumber
  \inttwo_\SDRDRD b(\omega,v_1 - v) \,
                  \domega \, M_1 \, \dv_1 \, M \, \dv = 1 \,.
\end{equation}
Because $\dmu=b(\omega,v_1-v)\,\domega\,M_1\dv_1\,M\dv$ is a positive
unit measure on $\SDRDRD$, we denote by $\lANGLE\Xi\rANGLE$ the
average over this measure of any integrable function
$\Xi=\Xi(\omega,v_1,v)$
\begin{equation}
  \label{ANGLE-Def}
  \lANGLE \Xi \rANGLE
  = \intthree_\SDRDRD \Xi(\omega,v_1,v) \, \dmu \,.
\end{equation}


\subsection{DiPerna-Lions-Mischler Theory}

   As in \cite{BGL, GL, LM}, we will work in the framework of
DiPerna-Lions solutions to the scaled Boltzmann equation on
the phase space $\OO=\RD\times\Omega$
\begin{equation}
  \label{Boltzmann-E}
\begin{aligned}
  \del_t G_\eps + v \DOT \GRAD G_\eps
  & = \frac{1}{\eps} \QQ(G_\eps,G_\eps) & \quad
  & \mbox{on} \quad \OO \times \R_+ \,,
\\
  \Ge(v,x,0) & = G_\eps^\init(v,x) & \quad
  & \mbox{on} \quad \OO \,,
\end{aligned}
\end{equation}
with the Maxwell reflection boundary condition \eqref{Boltzmann-BC}
which can be expressed as
\begin{equation}
  \label{Boltzmann-EBC}
  \gamma_{-} G_\eps 
  = (1 - \alpha) L(\gamma_{+} G_\eps)
    + \alpha \< \gamma_+ G_\eps \>_{\!\pO} \quad
  \mbox{on} \quad \Sigma_{-} \times \R_+ \,,
\end{equation}
where $\gamma_{\pm}G_\eps$ denote the traces of $G_\eps$ on the 
outgoing and incoming sets $\Sigma_\pm$.  Here the local 
reflection operator $L$ is defined to act on any 
$|v\cdot\nRm|M\dv\,\dsig$-measurable function $\phi$ over
$\del\OO$ by
\begin{equation}
  \nonumber
  L\phi(v,x) = \phi( \RRm(x) v,x) \quad
  \hbox{for almost every $(v,x)\in\del\OO$} \,,
\end{equation}
where $\RRm(x) v=v-2v\DOT\nRm(x)\nRm(x)$ is the specular 
reflection of $v$, while the diffuse reflection operator is 
defined as
\begin{equation}
  \nonumber
  \< \phi \>_{\pO} 
  = \sqrt{2\pi} \int_{v \cdot \nRm(x)>0}
                    \phi(v,x) \, v \DOT \nRm(x) \, M \dv \,.
\end{equation}
DiPerna-Lions theory requires that both the equation and boundary 
conditions in \eqref{Boltzmann-E} should be understood in the 
renormalized sense, see \eqref{scale-BE} and \eqref{renorm-bc}.  
These solutions were initially constructed by DiPerna and Lions 
\cite{DiPL} over the whole space $\RD$ for any initial data satisfying 
natural physical bounds.  For bounded domain case, Mischler \cite{Misch} 
recently developed a theory to treat the Maxwell reflection boundary 
condition \eqref{Boltzmann-EBC}.

   The DiPerna-Lions theory does not yield solutions that are known 
to solve the Boltzmann equation in the usual sense of weak solutions.  
Rather, it gives the existence of a global weak solution to a class 
of formally equivalent initial value problems that are obtained by 
multiplying \eqref{Boltzmann-E} by $\Gamma'(G_\eps)$, where $\Gamma'$ 
is the derivative of an admissible function $\Gamma$:
\begin{equation}
  \label{BE-Gamma}
  (\del_t + v \DOT \GRAD) \Gamma(G_\eps)
  = \frac{1}{\eps} \, \Gamma'(G_\eps) \QQ(G_\eps,G_\eps) \quad
  \mbox{on} \quad \OO\times\R_+ \,.
\end{equation}
Here a function $\Gamma:[0,\infty)\to\R$ is called admissible if it 
is continuously differentiable and for some $C_\Gamma<\infty$ its 
derivative satisfies
\begin{equation}
  \nonumber
  |\Gamma'(Z)| \leq \frac{C_\Gamma}{\sqrt{1 + Z}} \quad
  \hbox{for every $Z\in[0,\infty)$} \,.
\end{equation}
The solutions are nonnegative and lie in 
$C([0,\infty);\mbox{w-}L^1(M\dv\,\dx))$, where the prefix ``w-'' on
a space indicates that the space is endowed with its weak topology.

   Mischler \cite{Misch} extended DiPerna-Lions theory to domains 
with a boundary on which the Maxwell reflection boundary condition 
\eqref{Boltzmann-EBC} is imposed.  This required the proof of a
so-called trace theorem that shows that the restriction of $G_\eps$
to $\del\OO\times\R_+$, denoted $\gamma G_\eps$, makes sense.  In
particular, Mischler showed that $\gamma G_\eps$ lies in the set 
of all $|v\DOT\nRm|M\dv\,\dsig\,\dt$-measurable functions over 
$\del\OO\times\R_+$ that are finite almost everywhere, which we 
denote $L^0(|v\DOT\nRm|M\dv\,\dsig\,\dt)$.  He then defines 
$\gamma_\pm G_\eps= \One_{\Sigma_\pm} \gamma G_\eps$.  He proves 
the following.

\begin{Thm}
{\em (DiPerna-Lions-Mischler Renormalized Solutions \cite{DiPL, Misch})}  
Let $b$ be a collision kernel that satisfies the assumptions in 
Section \ref{AssumpKernel}.  Fix $\eps>0$.  Let $G_\eps^\init$ be  
any initial data in the entropy class
\begin{equation}
  \label{entropy-class}
  E(M\dv\,\dx) = \big\{ G_\eps^\init \geq 0 \,:\, 
                        H(G_\eps^\init) < \infty \big\} \,,
\end{equation} 
where the relative entropy functional is given by
\begin{equation}
  \nonumber
  H(G) = \int_\Omega \< \eta(G) \> \, \dx \quad \hbox{with} \quad
  \eta(G) = G \log(G) - G + 1 \,.
\end{equation}
Then there exists a $G_\eps\geq 0$ in 
$C([0,\infty);\mbox{w-}L^1(M\dv\,\dx))$ with $\gamma G_\eps\geq 0$
in $L^0(|v\cdot\nRm|M\dv\,\dsig\,\dt)$ such that:

\begin{itemize}

\item $G_\eps$ satisfies the global entropy inequality
\begin{equation}
  \label{entropy-inequality}
  H(G_\eps(t))
  + \int^t_0 \left[ \frac{1}{\eps} R(G_\eps(s))
                    + \frac{\alpha}{\sqrt{2\pi}}
                      \EE(\gamp G_\eps(s)) \right] \, \ds
  \leq H(G_\eps^\init) \quad
  \hbox{for every $t>0$} \,,
\end{equation}
where the entropy dissipation rate functional is given by
\begin{equation}
  \label{entropy-rate}
  R(G) = \frac{1}{4}
         \int_\Omega 
             \LANGLE \log\!\left( \frac{G'_1 G'}{G_1 G} \right)
                     \big( G'_1 G' - G_1 G \big) \RANGLE\, \dx \,,
\end{equation}
and the so-called Darroz\`{e}s-Guiraud information is given by
\begin{equation}
  \label{DG-inform}
  \EE(\gamp G)
  = \int_{\pO} 
        \left[ \big\< \eta(\gamp G) \big\>_{\pO}
               - \eta\big( \< \gamp G \>_{\pO} \big) 
               \right] \, \dsig \,;
\end{equation}
 
\item $G_\eps$ satisfies  
\begin{equation}
  \label{weak-BE-G}
\begin{gathered}
  \int_\Omega \< \Gamma( G_\eps(t_2)) \, Y \> \, \dx
  - \int_\Omega \< \Gamma( G_\eps(t_1)) \, Y \> \, \dx
  + \int^{t_2}_{t_1} \!\! \int_{\pO}
       \< \Gamma(\gamma G_\eps) \, Y \, (v \DOT \nRm) \> \, 
       \dsig \, \dt
\\
  - \int^{t_2}_{t_1} \!\! \int_\Omega
       \< \Gamma(G_\eps) \, v \DOT \GRAD Y \> \, \dx \, \dt
  = \frac{1}{\eps} \int^{t_2}_{t_1} \!\! \int_\Omega
       \< \Gamma'(G_\eps) \, 
          \mathcal{Q}(G_\eps,G_\eps) \, Y \> \, \dx \, \dt \,,
\end{gathered}
\end{equation}
for every admissible function $\Gamma$, every 
$Y\in C^1\cap L^\infty(\RD\times\bar{\Omega})$, and every
$[t_1,t_2]\subset[0,\infty]$;

\item $G_\eps$ satisfies
\begin{equation}
  \label{renormalized-BC}
   \gamn G_\eps = (1 - \alpha) L(\gamp G_\eps)
                  + \alpha \< \gamp G_\eps \>_{\pO} \quad
  \mbox{almost everywhere on} \quad  \Sigma_- \times \R_+ \,.
\end{equation}
\end{itemize}
\end{Thm}

\smallskip
\noindent
{\bf Remark.} Because the trace $\gamma G_\eps$ is only known to 
exist in $L^0(|v\cdot\nRm|M\dv\,\dsig\,\dt)$ rather than in
$L^1_{loc}(\dt;L^1(|v\DOT\nRm|M\dv\,\dsig))$, we cannot conclude
from the boundary condition \eqref{renormalized-BC} that 
\begin{equation}
  \label{normal-vanish}
  \< v \, \gamma G_\eps \> \DOT \nRm = 0 \quad 
  \mbox{on} \quad \del\Omega \,.
\end{equation}
Indeed, we cannot even conclude that the boundary mass-flux
$\< v\,\gamma G_\eps\>\DOT\nRm$ is defined on $\del\Omega$.  
Moreover, in contrast to DiPerna-Lions theory over the whole
space or periodic domains, it is not asserted in \cite{Misch} 
that $G_\eps$ satisfies the weak form of the local mass 
conservation law
\begin{equation}
  \label{weak-mass-G}
  \int_\Omega \chi \, \< G_\eps(t_2) \> \, \dx
  - \int_\Omega \chi \, \< G_\eps(t_1) \> \, \dx
  - \int^{t_2}_{t_1} \!\! \int_\Omega
       \GRAD \chi \DOT \< v \, G_\eps \> \, \dx \, \dt
  = 0 \qquad \forall \chi \in C^1(\overline{\Omega}) \,.
\end{equation}
If this were the case, it would allow a great simplification
the proof of our main result.  Rather, we will employ the boundary 
condition \eqref{renormalized-BC} inside an approximation to
\eqref{weak-mass-G} that has a well-defined boundary flux.

\smallskip
\noindent
{\bf Remark.}  As was shown in \cite{BGL3}, the condition 
$H(G) < \infty$ found in our definition of the entropy class 
\eqref{entropy-class} is equivalent to the conditon
$$
  \inttwo_{\RD\times\Omega} 
      \big( 1 + |v|^2 + |\log(G)| \big) \, 
      G \, M \dv \, \dx < \infty \,, 
$$
which is used by Mischler and others.  By presenting it as we do, it 
is clear that the entropy class is simply those kinetic densities $G$ 
whose relative entropy with respect to $M$ is finite. 


\section{Main Result}
\label{MainResult}

\subsection{Main Theorem}

   We will consider families $G_\eps$ of DiPerna-Lions renormalized
solutions to \eqref{Boltzmann-E} such that $G_\eps^\init\geq 0$
satisfies the entropy bound
\begin{equation}
  \label{entropy-bound}
  H(G_\eps^\init) \leq C^\init \delta_\eps^{\,2}
\end{equation}
for some $C^\init<\infty$ and $\delta_\eps>0$ that satisfies the
scaling $\delta_\eps\to 0$ as $\eps\to 0$.

   The value of $H(G)$ provides a natural measure of the proximity of 
$G$ to the equilibrium $G=1$.  We define the families $g^\ini_\eps$ 
and $g_\eps$ of fluctuations about $G=1$ by the relations
\begin{equation}
  \label{fluctuations}
  G^\ini_\eps = 1 + \dep g^\ini_\eps \,, \qquad
  G_\eps = 1 + \dep g_\eps \,.
\end{equation}
One easily sees \cite{BGL3} that $H$ asymptotically behaves like half 
the square of the $L^2\mbox{-}$norm of these fluctuations as $\eps\to0$.
Hence, the entropy bound \eqref{entropy-bound} combined with the
entropy inequality \eqref{entropy-inequality} is consistent with 
these fluctuations being of order $1$.  Just as the relative entropy $H$ 
controls the fluctuations $g_\eps$, the dissipation rate $R$ given by
\eqref{entropy-rate} controls the scaled collision integrals defined by
\begin{equation}
  \nonumber
  q_\eps = \frac{1}{\sqrt{\eps} \delta_\eps}
           \big( G'_{\eps 1} G'_\eps - G_{\eps 1} G_\eps \big) \,.
\end{equation}

   Here we only state the weak acoustic limit theorem because the
corresponding strong limit theorem is analogous to that stated in 
\cite{GL} and its proof based on the weak limit theorem and relative 
entropy convergence is essentially the same.

\begin{Thm}
  \label{weak-acoustic}
{\em (Weak Acoustic Limit Theorem)} Let $b$ be a collision kernel
that satisfies the assumptions in Section \ref{AssumpKernel}.
Let $G_\eps^\init$ be a family in the entropy class $E(M\dv\,\dx)$ 
that satisfies the normalization \eqref{norm-init} and the entropy 
bound \eqref{entropy-bound} for some $C^\init<\infty$ and 
$\delta_\eps>0$ satisfies the scaling
\begin{equation}
  \nonumber
  \delta_\eps = O\big( \sqrt{\eps} \big) \,.
\end{equation}
Assume, moreover, that for some
$(\rho^\init,u^\init,\theta^\init)\in L^2(\dx;\RRDR)$
the family of fluctuations $g_\eps^\init$ defined by
\eqref{fluctuations} satisfies
\begin{equation}
  \label{initial-limit}
  (\rho^\init, u^\init, \theta^\init)
  = \lim_{\eps\to 0}
        \left( \< g_\eps^\init \> \,,\,
               \< v \, g_\eps^\init \> \,,\,
               \big\< \big( \tfrac{1}{\D} |v|^2 - 1 \big) \,
                            g_\eps^\init \big\> \right) \quad
  \mbox{in the sense of distributions} \,.
\end{equation}

   Let $G_\eps$ be any family of DiPerna-Lions-Mischler 
renormalized solutions to the Boltzmann equation 
\eqref{Boltzmann-E} that have $G^\init_\eps$ as initial values.

   Then, as $\eps\to 0$, the family of fluctuations $g_\eps$ defined 
by \eqref{fluctuations} satisfies
\begin{equation}
  \label{lim-fluct}
  g_\eps \to \rho + v \DOT u
                  + (\tfrac{1}{2} |v|^2 - \tfrac{\D}{2}) \theta \quad
  \hbox{in $\wL^1_{loc}(\dt;\wL^1((1+|v|^2)M\dv\,\dx))$} \,,
\end{equation}
where $(\rho,u,\theta)\in C([0,\infty);L^2(\dx;\RRDR))$ is the
unique solution to the acoustic system \eqref{Acoustic-IVP} that
satisfies the impermeable boundary condition \eqref{Acoustic-BC} 
and has initial data $(\rho^\init,u^\init,\theta^\init)$ obtain from 
\eqref{initial-limit}.
In addition, $\rho$ satisfies
\begin{equation}
  \label{lim-mass-tot}
  \int_\Omega \rho \, \dx = 0 \,.
\end{equation}
\end{Thm}

   This result improves upon the acoustic limit result in \cite{GL}
in three ways.  First, its assumption on the collision kernel $b$ is
the same as \cite{LM}, so it treats a broader class of cut-off kernels
than was treated in \cite{GL}.  In particular, it treats kernals 
derived from soft potentials.  Second, its scaling assumption is
$\dep=O(\sqrt{\eps})$, which is certainly better than the scaling 
assumption (\ref{GL-scaling}) used in \cite{GL}.  This assumption is
still a long way from that required by the formal derivation of the
acoustic system, which is $\dep\to 0$ as $\eps\to 0$.  Our more 
restrictive requirement arises from the way in which we remove the 
local conservation law defects of the DiPerna-Lions solutions.  
Third, we derive a weak form of the boundary condition $u\DOT\nRm=0$.  
It is the first time such a boundary condition for the acoustic system 
is derived from the Boltzmann equation with the Maxwell reflection 
boundary condition.


\subsection{Proof of the Main Theorem}

   In order to derive the fluid equations with boundary conditions,
we need to pass to the limit in approximate local conservation laws
built from the renormalized Boltzmann equation \eqref{BE-Gamma}.
We choose the renormalization used in \cite{LM} --- namely,
\begin{equation}
  \label{renorm}
  \Gamma(Z) = \frac{Z-1}{1+(Z-1)^2} \,.
\end{equation}
After dividing by $\dep$, equation \eqref{BE-Gamma} becomes
\begin{equation}
  \label{scale-BE}
  \del_t \gps + v \DOT \GRAD \gps
  = \frac{1}{\sqrt{\eps}} \, \Gamma'(G_\eps)
    \inttwo_\SDRD  q_\eps \, b(\omega\,,v_1-v) \, 
                   \domega \, M_1 \dv_1 \,,
\end{equation}
where $\gps=\Gamma(G_\eps)/\dep$.  By introducing
$N_\eps=1+\dep^2g^2_\eps$, we can write
\begin{equation}\label{renorm-g}
  \gps = \frac{g_\eps}{N_\eps} \,, \qquad
  \Gamma'(G_\eps) = \frac{2}{N^2_\eps} - \frac{1}{N_\eps} \,.
\end{equation}

   When moment of the renormalized Boltzmann equation
\eqref{scale-BE} is formally taken with respect to any
$\zeta\in\mbox{span}\{1\,,v_1\,,\cdots\,,v_\D\,,|v|^2\}$, one obtains
\begin{equation}
  \label{moment-equation}
  \del_t \< \zeta \, \gps \>
  + \DIV \< v \, \zeta \, \gps \>
  = \frac{1}{\sqrt{\eps}} \,
    \lANGLE \zeta \, \Gamma'(G_\eps) \, q_\eps \rANGLE \,.
\end{equation}
This fails to be a local conservation law because the so-called
{\em conservation defect} on the right-hand side is generally nonzero.
We will show that this defect vanishes as $\eps\to 0$, while the 
left-hand side converges to the local conservation law corresponding 
to $\zeta$.  More precisely, it can be shown that every DiPerna-Lions 
solution satisfies \eqref{moment-equation} in the sense that for 
every $\chi\in C^1(\Omega)$ and every $[t_1,t_2]\subset[0,\infty)$ it
satisfies
\begin{equation}
  \label{weak-moment}
\begin{aligned}
  \int_\Omega \chi \, \< \zeta \, \gps(t_2) \> \, \dx
  & - \int_\Omega \chi \, \< \zeta \, \gps(t_1) \> \, \dx
    + \int_{t_1}^{t_2} \!\! \int_{\del\Omega}
          \chi \, \< v \, \zeta \, \gamma \gps \> \DOT \nRm \,
          \dsig \, \dt
\\
  & - \int_{t_1}^{t_2} \!\! \int_\Omega
          \GRAD \chi \DOT \< v \, \zeta \, \gps \> \, \dx \, \dt
    = \int_{t_1}^{t_2} \!\! \int_\Omega
          \chi \, \frac{1}{\sqrt{\eps}} \,
          \lANGLE \zeta \, \Gamma'(G_\eps) \, q_\eps \rANGLE \,
          \dx \, \dt \,.
\end{aligned}
\end{equation}
Moreover, from \eqref{renormalized-BC} the boundary condition 
is understood in the renormalized sense:
\begin{equation}
  \label{renorm-bc}
   \gamn \gps
   =  \frac{(1 - \alpha) L \gamp g_\eps 
            + \alpha \< \gamp g_\eps \>_{\pO}}
           {1 + \delta^2_\eps [(1 - \alpha) L \gamp g_\eps 
                               + \alpha \< \gamp g_\eps\>_{\pO}]^2} 
   \quad \mbox{on $\Sigma_- \times \R_+$} \,,
\end{equation}
where the equality holds almost everywhere.  We will pass to the limit
in the weak form \eqref{weak-moment}.  The Main Theorem will be proved 
in two steps: the interior equations will be established first and the 
boundary condition second.  

   The acoustic system \eqref{Acoustic-IVP} is justified in the 
interior of $\Omega$ by showing that the limit of \eqref{weak-moment} 
as $\eps\to0$ is the weak form of the acoustic system whenever the 
test function $\chi$ vanishes on $\del\Omega$.   We prove that the 
conservation defect on the right-hand side of \eqref{weak-moment} 
vanishes as $\eps\to0$ in Proposition \ref{Defect-Prop}, which is 
presented in the next section.  The proof of the analogous result 
in \cite{GL} must be modified in order to include the case 
$\dep=O(\sqrt{\eps})$.  The convergence of the density and flux 
terms is proved essentially as in \cite{GL}, so we omit those 
arguments here.  The upshot is that every converging subsequence of 
the family of fluctuations $g_\eps$ satisfies
\begin{equation}
  \nonumber
  g_\eps \to \rho + v \DOT u
                  + (\tfrac{1}{2} |v|^2 - \tfrac{\D}{2}) \theta \quad
  \hbox{in $\wL^1_{loc}(\dt;\wL^1((1+|v|^2)M\dv\,\dx))$} \,,
\end{equation}
where $(\rho,u,\theta)\in C([0,\infty);\wL^2(\dx;\RRDR))$ satisfies
for every $[t_1,t_2]\subset[0,\infty)$ 
\begin{subequations}
\begin{align}
  \label{lim-mass}
  \int_\Omega \chi \, \rho(t_2) \, \dx 
  - \int_\Omega \chi \, \rho(t_1) \, \dx 
  - \int_{t_1}^{t_2} \!\!
        \int_\Omega \GRAD \chi \DOT u \, \dx \, \dt 
  & = 0 \qquad \hbox{$\forall \chi\in C_0^1(\overline{\Omega})$} \,,
\\
  \label{lim-motion}
  \int_\Omega w \DOT u(t_2) \, \dx 
  - \int_\Omega w \DOT u(t_1) \, \dx 
  - \int_{t_1}^{t_2} \!\! 
        \int_\Omega \DIV w \, (\rho + \theta) \, \dx \, \dt 
  & = 0 \qquad \hbox{$\forall w\in C_0^1(\overline{\Omega};\RD)$} \,,
\\
  \label{lim-heat}
  \tfrac{\D}2 \int_\Omega \chi \, \theta(t_2) \, \dx 
  - \tfrac{\D}2 \int_\Omega \chi \, \theta(t_1) \, \dx 
  - \int_{t_1}^{t_2} \!\!
        \int_\Omega \GRAD \chi \DOT u \, \dx \, \dt 
  & = 0 \qquad \hbox{$\forall \chi\in C_0^1(\overline{\Omega})$} \,.
\end{align}
\end{subequations}
This shows that the acoustic system \eqref{Acoustic-IVP} is satisfied
in the interior of $\Omega$.

   The more significant step is to justify the impermeable boundary 
condition \eqref{Acoustic-BC}.  Unlike to what is done for the 
incompressible Stokes \cite{MS} and Navier-Stokes \cite{JM} limits, 
here we do not have enough control to pass to the limit in the 
boundary terms in \eqref{weak-moment} for the local conservation 
laws of momentum and energy.  We can however do so for the local 
conservation law of mass --- i.e. when $\zeta=1$.  Indeed, 
Proposition \ref{Boundary-Prop} of Section 5 will show that we 
can extend \eqref{lim-mass} to
\begin{equation}
  \label{lim-mass-bdry}
  \int_\Omega \chi \, \rho(t_2) \, \dx 
  - \int_\Omega \chi \, \rho(t_1) \, \dx 
  - \int_{t_1}^{t_2} \!\!
        \int_\Omega \GRAD \chi \DOT u \, \dx \, \dt 
  = 0 \qquad \hbox{$\forall \chi\in C^1(\overline{\Omega})$} \,.
\end{equation}
We obtain \ref{lim-mass-tot} by setting $\chi=1$ and $t_1=0$ above,
and using the fact that the family $G_\eps^\init$ satisfies the
normalization \eqref{norm-init}.  

   Because for every $\chi\in C^1(\overline{\Omega})$ we can find a 
sequence $\{\chi_n\}\subset C^1_0(\overline{\Omega})$ such that 
$\chi_n\to\chi$ in $L^2(\dx)$, it follows from \eqref{lim-mass} and 
\eqref{lim-heat} that 
\begin{equation}
  \nonumber
\begin{aligned}
  \tfrac{\D}2 \int_\Omega \chi \, \theta(t_2) \, \dx 
  - \tfrac{\D}2 \int_\Omega \chi \, \theta(t_1) \, \dx 
  & = \lim_{n\to\infty} 
          \tfrac{\D}2 \int_\Omega \chi_n \, \theta(t_2) \, \dx 
      - \lim_{n\to\infty} 
            \tfrac{\D}2 \int_\Omega \chi_n \, \theta(t_1) \, \dx 
\\
  & = \lim_{n\to\infty} 
          \int_\Omega \chi_n \, \rho(t_2) \, \dx 
      - \lim_{n\to\infty} 
            \int_\Omega \chi_n \, \rho(t_1) \, \dx 
\\
  & = \int_\Omega \chi \, \rho(t_2) \, \dx 
      - \int_\Omega \chi \, \rho(t_1) \, \dx \,.
\end{aligned}
\end{equation}
It thereby follows from \eqref{lim-mass-bdry} that we can extend 
\eqref{lim-heat} to
\begin{equation}
  \label{lim-heat-bdry}
  \tfrac{\D}2 \int_\Omega \chi \, \theta(t_2) \, \dx 
  - \tfrac{\D}2 \int_\Omega \chi \, \theta(t_1) \, \dx 
  - \int_{t_1}^{t_2} \!\!
        \int_\Omega \GRAD \chi \DOT u \, \dx \, \dt 
  = 0 \qquad \hbox{$\forall \chi\in C^1(\overline{\Omega})$} \,.
\end{equation}

   Finally, because for every $w\in C^1(\overline{\Omega};\RD)$ 
such that $w\DOT\nRm=0$ on $\del\Omega$ we can find a sequence 
$\{w_n\}\subset C^1_0(\overline{\Omega};\RD)$ such that 
$w_n\to w$ in $L^2(\dx;\RD)$ and $\DIV w_n\to\DIV w$ in 
$L^2(\dx)$, it follows from \eqref{lim-motion} that
\begin{equation}
  \nonumber
\begin{aligned}
  \int_\Omega w \DOT u(t_2) \, \dx 
  - \int_\Omega w \DOT u(t_1) \, \dx 
  & = \lim_{n\to\infty} \int_\Omega w_n \DOT u(t_2) \, \dx 
      - \lim_{n\to\infty} \int_\Omega w_n \DOT u(t_1) \, \dx
\\
  & = \lim_{n\to\infty}
          \int_{t_1}^{t_2} \!\! 
          \int_\Omega \DIV w_n \, (\rho + \theta) \, \dx \, \dt 
\\
  & = \int_{t_1}^{t_2} \!\! 
          \int_\Omega \DIV w \, (\rho + \theta) \, \dx \, \dt \,.
\end{aligned}
\end{equation}
But this combined with \eqref{lim-mass-bdry} and \eqref{lim-heat-bdry}
is the weak formulation of the acoustic system \eqref{Acoustic-IVP}
with the boundary condition \eqref{Acoustic-BC}.  Because this system
has a unique weak solution in $C([0,\infty);\wL^2(\dx;\RRDR))$, all
converging sequences of the family $g_\eps$ have this same limit.  
Moreover, this limit must be the strong solution that lies in 
$C([0,\infty);L^2(\dx;\RRDR))$.  The family of fluctuations $g_\eps$ 
therefore converges as asserted by \eqref{lim-fluct}.  \qed

\smallskip
\noindent
{\bf Remark.} Had we known that $G_\eps$ satisfies the weak form 
of the local mass conservation law \eqref{weak-mass-G} then we 
could have easily obtained \eqref{lim-mass-bdry} by passing to the 
limit in \eqref{weak-mass-G}.  In that case there would be a great 
simplification in our proof because there would be no need for 
Proposition \ref{Boundary-Prop}.  


\section{Removal of the Conservation Defects}

   The conservation defects in \eqref{moment-equation} have the form 
\begin{equation}
  \nonumber
  \frac{1}{\sqrt{\eps}} \,
  \lANGLE \zeta \, \Gamma'(G_\eps) \, q_\eps \rANGLE
  = \frac{1}{\sqrt{\eps}} \,
    \LANGLE \zeta \,
            \bigg( \frac2{N_\eps^{\,2}} - \frac1{N_\eps} \bigg) \,
            q_\eps \RANGLE \,.
\end{equation}
In order to establish local conservation laws, we must show that 
these defects vanish as $\eps\to0$.  This is done with the following 
proposition.

\begin{Prop}
  \label{Defect-Prop}
For $n=1$ and $n=2$, and for every
$\zeta\in \mbox{span}\{ 1, v_1, \cdots, v_{\D}, |v|^2\}$ one has
\begin{equation}
  \label{vanishing}
  \frac{1}{\sqrt{\eps}} \,
  \LANGLE \zeta \, \frac{q_\eps}{N_\eps^{\,n}} \RANGLE \to 0 \quad
  \hbox{in $\wL^1_{loc}(\dt;\wL^1(\dx))$ as $\eps\to 0$} \,.
\end{equation}
\end{Prop}

\begin{proof}
Similar to the proof of Proposition 8.1 in \cite{LM}, for $n=1$,
we obtain the decomposition
\begin{equation}
  \label{decomp-1}
  \frac{1}{\sqrt{\eps}}
  \LANGLE \zeta \, \frac{q_\eps}{N_\eps} \RANGLE
  = \frac{\dep^2}{\sqrt{\eps}} \,
    \LANGLE \zeta \, \frac{g_{\eps 1}^{\,2} q_\eps}
                          {N_{\eps 1}N_\eps} \RANGLE
    + \LANGLE \zeta \, 
              \frac{\delta_\eps^{\,2} 
                    (g_{\eps 1} + g_\eps) q_\eps^{\,2}}
                   {N'_{\eps 1} N'_\eps N_{\eps 1} N_\eps} \RANGLE
    - \frac{\delta_\eps^{\,2}}{\sqrt{\eps}} \,
      \LANGLE \zeta' \, \frac{g'_{\eps 1} g'_\eps \, q_\eps}
                             {N'_{\eps 1} N'_\eps N_{\eps 1} N_\eps} \,
                        J_\eps \RANGLE \,,
\end{equation}
where $J_\eps$ is given by
\begin{equation}
  \label{J-Def}
  J_\eps = 2 + \dep \big( g'_{\eps 1} + g'_\eps 
                          + g_{\eps 1} + g_\eps \big)
             - \dep^2 \big( g'_{\eps 1} g'_\eps 
                            - g_{\eps 1} g_\eps \big) \,.
\end{equation}

   We can then dominate the integrands of the three terms on the
right-hand side of \eqref{decomp-1}.  Because for every
$\zeta\in \mbox{span}\{ 1, v_1, \cdots, v_{\D}, |v|^2\}$
there exists a constant $C<\infty$ such that $|\zeta| \leq C \sigma$
where $\sigma \equiv1+|v|^2$, the integrand of the first term is dominated by
\begin{equation}
  \label{term-1}
  \frac{\dep^2}{\sqrt{\eps}} \,
  \sigma \, \frac{g^2_{\eps 1} |q_\eps|}{N_{\eps 1} N_\eps} \,.
\end{equation}
Because
$\frac{\dep|g_{\eps 1}+ g_\eps|}{\sqrt{N'_{\eps 1}N'_\eps N_{\eps 1}N_\eps}} \leq 2$,
the integrand of the second term is dominated by
\begin{equation}
  \label{term-2}
  \sigma \,
  \frac{\dep q^2_\eps}{\sqrt{N'_{\eps 1} N'_\eps N_{\eps 1} N_\eps}} \,.
\end{equation}
Finally, because
$\frac{|J_\eps|}{\sqrt{N'_{\eps 1}N'_\eps N_{\eps 1}N_\eps}}\leq 8$,
the integrand of the third term is dominated by
\begin{equation}
  \label{term-3}
  \frac{\dep^2}{\sqrt{\eps}} \, \sigma' \,
  \frac{|g'_{\eps 1} g'_\eps| |q_\eps|}
       {\sqrt{N'_{\eps 1} N'_\eps N_{\eps 1} N_\eps}} \,.
\end{equation}
Hence, the result \eqref{vanishing} for the case $n=1$ will follow
once we establish that the terms \eqref{term-1}, \eqref{term-2}, and
\eqref{term-3} vanish as $\eps\to 0$.

   The result \eqref{vanishing} for the case $n=2$ will follow 
similarly.  We start with the decomposition
\begin{equation}
  \nonumber
\begin{aligned}
  \frac{1}{\sqrt{\eps}}
  \LANGLE \zeta \, \frac{q_\eps}{N_\eps^{\,2}} \RANGLE
  & = \frac{\dep^2}{\sqrt{\eps}} \,
      \LANGLE \zeta \, \frac{g_{\eps 1}^{\,2} q_\eps}
                            {N_{\eps 1}N_\eps}
              \bigg( 1 + \frac{1}{N_{\eps 1}} \bigg) \RANGLE
    + \LANGLE \zeta \, 
              \frac{\delta_\eps^{\,2} 
                    (g_{\eps 1} + g_\eps) q_\eps^{\,2}}
                   {N'_{\eps 1} N'_\eps N_{\eps 1} N_\eps}
              \bigg( \frac{1}{N'_{\eps 1} N'_\eps}
                     + \frac{1}{N_{\eps 1} N_\eps} \bigg) \RANGLE
\\
  & \quad \,
    - \frac{\delta_\eps^{\,2}}{\sqrt{\eps}} \,
      \LANGLE \zeta' \, \frac{g'_{\eps 1} g'_\eps \, q_\eps}
                             {N'_{\eps 1} N'_\eps N_{\eps 1} N_\eps} \,
                        J_\eps
              \bigg( \frac{1}{N'_{\eps 1} N'_\eps}
                     + \frac{1}{N_{\eps 1} N_\eps} \bigg) \RANGLE \,,
\end{aligned}
\end{equation}
where $J_\eps$ is given by \eqref{J-Def}.  Because the terms in 
parentheses above are each bounded by $2$, we can dominate the 
three terms on the right-hand side above just as we did the terms 
on the right-hand side of \eqref{decomp-1} for the case $n=1$.  
The result \eqref{vanishing} for the case $n=2$ will then also 
follow once we establish that the terms \eqref{term-1}, 
\eqref{term-2}, and \eqref{term-3} vanish as $\eps\to 0$.

   That term \eqref{term-2} vanishes is easy to see.  The 
inequality
$n'_{\eps 1}n'_\eps n_{\eps 1}n_\eps \leq 2\sqrt{N'_{\eps 1}N'_\eps N_{\eps 1}N_\eps}$,
where $n_\eps=1+\frac{\dep}{3}g_\eps$, along with the estimate
\begin{equation}
  \nonumber
  \sigma \,
  \frac{q_\eps^{\,2}}{n'_{\eps 1} n'_\eps n_{\eps 1} n_\eps}
  = O\left( \left| \log\left( \sqrt{\eps} \delta_\eps \right) \right| \right) 
  \quad \mbox{in $L^1_{loc}(\dt;L^1(\dmu\,\dx))$ as $\eps\to0$} \,,
\end{equation}
which is proved in Lemma 9.4 of \cite{GL}, imply that
\begin{equation}
  \nonumber
  \sigma \, \frac{\delta_\eps q^2_\eps}
                 {\sqrt{N'_{\eps 1} N'_\eps N_{\eps 1} N_\eps}} 
  = O\left( \delta_\eps 
            \left| \log\left( \sqrt{\eps} \delta_\eps \right) \right| \right) 
  \to 0 \quad 
  \mbox{in $L^1_{loc}(\dt;L^1(\dmu\,\dx))$ as $\eps\to0$} \,.
\end{equation}
The fact that the terms \eqref{term-1} and \eqref{term-3} vanish as 
$\eps\to0$ follows from Lemma \ref{vanish-1}, which is proved below.  
We thereby complete the proof of Proposition \ref{Defect-Prop}.
\end{proof}

\begin{Lem}
  \label{vanish-1}
\begin{align}
  \label{term1}
  \frac{\delta_\eps^2}{\sqrt{\eps}} \, \sigma \,
  \frac{g^2_{\eps 1} |q_\eps|}{N_{\eps 1} N_\eps} \to 0 \quad
  & \mbox{in $L^1_{loc}(\dt;L^1(\dmu\,\dx))$ as $\eps\to0$} \,,
\\
  \label{term3}
  \frac{\dep^2}{\sqrt{\eps}} \, \sigma' \,
  \frac{|g'_{\eps 1} g'_\eps| |q_\eps|}
       {\sqrt{N'_{\eps 1} N'_\eps N_{\eps 1} N_\eps}} \to 0 \quad
  & \mbox{in $L^1_{loc}(\dt;L^1(\dmu\,\dx))$ as $\eps\to0$} \,.
\end{align}
\end{Lem}

\begin{proof}
The key to proving Lemma \ref{vanish-1} is the fact that
\begin{equation}
  \label{compactness}
  \frac{g_\eps^{\,2}}{\sqrt{N_\eps}} \quad
  \mbox{is relatively compact in} \!\!\quad
  w\mbox{-}L^1_{loc}(\dt; w\mbox{-}L^1(aM\dv\,\dx)) \,.
\end{equation}
This fact follows from Proposition 7.1 of \cite{LM}, where it 
plays an essential role in establishing the Navier-Stokes limit.
The approach to proving Lemma \ref{vanish-1} is the same used to
prove the analogous result in \cite{GL}.  There the terms 
\eqref{term-1} and \eqref{term-3} were estimated by using the 
entropy dissipation bound along with the nonlinear estimate
\begin{equation}
  \nonumber
  \sigma \, \frac{g_\eps^{\,2}}{\sqrt{N_\eps}} 
  = O(|\log(\delta_\eps)|) \quad
  \mbox{in $L^\infty(\dt;L^1(M\dv\,\dx))$ as $\eps \to 0$} \,.
\end{equation}
Here this nonlinear estimate, which originated in \cite{BGL3}, is 
replaced by the new weak compactness result \eqref{compactness}
from \cite{LM}, thereby extending the result in \cite{GL} to 
the scaling $\dep=\sqrt{\eps}$.

   The entropy inequality \eqref{entropy-inequality} and the 
entropy bound \eqref{entropy-bound} combine to bound the entropy 
dissipation as
\begin{equation}
  \label{dissipation-Bnd}
  \frac{1}{\eps \delta_\eps^{\,2}} \int^\infty_0 \!\! \int_\Omega 
  \LANGLE \frac{1}{4} 
          r\!\left( \frac{\sqrt{\eps} \delta_\eps q_\eps}
                         {G_{\eps 1} G_\eps} \right)
          G_{\eps 1} G_\eps \RANGLE \, \dx \, \dt
  \leq C^\init \,,
\end{equation}
where the function $r$ is defined over $z>-1$ by $r(z)=z\log(1+z)$.
The function $r$ is strictly convex over $z>-1$.  The proofs of
\eqref{term1} and \eqref{term3} are each based on a delicate use of
the classical Young inequality satisfied by $r$ and its Legendre
dual $r^*$, namely, the inequality
\begin{equation}
  \nonumber
  p z \leq r^*(p) + r(z) \quad
  \mbox{for every $p \in \R$ and $z> -1$} \,.
\end{equation}
For every positive $\varrho$ and $y$ we set
\begin{equation}
  \nonumber
  p = \frac{\sqrt{\eps} \dep y}{\varrho} \quad \mbox{and} \quad
  z = \frac{\sqrt{\eps} \dep |q_\eps|}{G_{\eps 1} G_\eps} \,,
\end{equation}
and use the fact that $r(|z|)\leq r(z)$ for every $z> -1$ to obtain
\begin{equation}
  \label{Young}
  y |q_\eps|
  \leq \frac{\varrho}{\eps \delta_\eps^{\,2}} \,
       r^*\!\left( \frac{\sqrt{\eps} \delta_\eps y}{\varrho} \right) 
       G_{\eps 1} G_\eps
       + \frac{\varrho}{\eps \delta_\eps^{\,2}} \,
         r\!\left( \frac{\sqrt{\eps} \delta_\eps q_\eps}
                        {G_{\eps 1} G_\eps} \right) 
         G_{\eps 1} G_\eps \,.
\end{equation}
This inequality is the starting point for the proofs of assertions
\eqref{term1} and \eqref{term3}.  These proofs also use the facts, 
recalled from \cite{BGL3}, that $r^*$ is superquadratic in the sense
\begin{equation}
  \label{superquadratic}
  r^*(\lambda p) \leq \lambda^2 r^*(p) \quad
  \mbox{for every $p>0$ and $\lambda\in[0,1]$} \,,
\end{equation}
and that $r^*$ has the exponential asymptotics $r^*(p)\sim\exp(p)$ as
$p\to\infty$.

   The proof of assertion \eqref{term1} follows that of Lemma 8.2 
in \cite{LM}.  We use the inequality \eqref{Young} with 
$y=\frac{\sigma}{4s^*}\frac{\delta_\eps^{\,2}}
{\sqrt{\eps}}\frac{g_{\eps 1}^{\,2}}{N_{\eps 1}N_\eps}$, where 
$s^*\in[1,\infty)$ is related to $s\in(1,\infty]$ appearing in
\eqref{bBnd} by the duality relation $\frac1s+\frac1{s^*}=1$.
We then apply the superquadratic property \eqref{superquadratic}
with $\lambda=\frac{\dep^3 g^2_{\eps 1}}{\varrho N_{\eps 1}N_\eps}$
and $p=\frac{\sigma}{4 s^*}$, where we note that $\lambda \leq 1$
whenever $\delta_\eps\leq\varrho$.  This leads to the bound
\begin{equation}
  \label{term1-Bnd}
  \frac{\sigma}{4 s^*} \, 
  \frac{\delta_\eps^{\,2}}{\sqrt{\eps}} \,
  \frac{g^2_{\eps 1} |q_\eps|}{N_{\eps 1} N_\eps} 
  \leq \frac{1}{\varrho} \, \frac{\delta_\eps^{\,4}}{\eps} \,
       \frac{g_{\eps 1}^{\,4}}{N_{\eps 1}^{\,2} N_\eps^{\,2}} \,
       r^*\!\left( \frac{\sigma}{4 s^*} \right)
       G_{\eps 1} G_\eps
       + \frac{\varrho}{\eps \delta_\eps^{\,2}} 
         r\!\left( \frac{\sqrt{\eps} \delta_\eps q_\eps}
                        {G_{\eps 1} G_\eps} \right) 
         G_{\eps 1} G_\eps \,.
\end{equation}
The second term on the right-hand side above can be made 
arbitrarily small in $L^1(\dmu\,\dx\,\dt)$ by using the
entropy dissipation bound \eqref{dissipation-Bnd} and 
picking $\varrho$ small enough.  Assertion \eqref{term1} will 
then follow upon showing that for every $\varrho>0$ the first 
term on the right-hand side of \eqref{term1-Bnd} vanishes as 
$\eps\to0$. 

   Because $G_{\eps 1}G_\eps\leq 2\sqrt{N_{\eps 1}N_\eps}$ while
$N_\eps\geq 1$, the first term on the right-hand side of 
\eqref{term1-Bnd} is bounded by
\begin{equation}
  \nonumber
  \frac{2 \, \delta_\eps^{\,2}}{\varrho \, \eps} \,
  \frac{\delta_\eps^{\,2} g_{\eps 1}^{\,2}}{N_{\eps 1}} \,
  \frac{g_{\eps 1}^{\,2}}{\sqrt{N_{\eps 1}}} \,
  r^*\!\left( \frac{\sigma}{4 s^*} \right) \,.
\end{equation}
The first factor above is bounded because 
$\delta_\eps^{\,2}=O(\eps)$, while the second is bounded
above by $1$ and satisfies
\begin{equation}
  \nonumber
  \frac{\delta_\eps^{\,2} g_{\eps 1}^{\,2}}{N_{\eps 1}} 
  \to 0 \quad 
  \hbox{in measure as $\eps\to 0$} \,.  
\end{equation}
It follows from \eqref{compactness} and Lemma 8.1 of \cite{LM}
that 
\begin{equation}
  \nonumber
  \frac{g_{\eps 1}^{\,2}}{\sqrt{N_{\eps 1}}} \,
  r^*\!\left( \frac{\sigma}{4 s^*} \right) \quad
  \hbox{is relatively compact in} \quad
  w\mbox{-}L^1_{loc}(\dt; w\mbox{-}L^1(\dmu\,\dx)) \,.
\end{equation}
We thereby conclude by the Product Limit Theorem \cite{BGL3} that
\begin{equation}
  \nonumber
  \frac{2 \, \delta_\eps^{\,2}}{\varrho \, \eps} \,
  \frac{\delta_\eps^{\,2} g_{\eps 1}^{\,2}}{N_{\eps 1}} \,
  \frac{g_{\eps 1}^{\,2}}{\sqrt{N_{\eps 1}}} \,
  r^*\!\left( \frac{\sigma}{4 s^*} \right) \to 0 \quad
  \mbox{in $L^1_{loc}(\dt;L^1(\dmu\,\dx))$} \,. 
\end{equation}
Hence, for every $\varrho>0$ the first term on the right-hand 
side of \eqref{term1-Bnd} vanishes as $\eps\to0$.  Assertion
\eqref{term1} thereby follows.

   The proof of assertion \eqref{term3} similarly follows that 
of Lemma 8.3 in \cite{LM}.  We use the inequality \eqref{Young}
with $y=\frac{\sigma}{4 s^*}\frac{\dep^2}{\sqrt{\eps}}
\frac{|g'_{\eps 1}g'_\eps|}{\sqrt{N'_{\eps 1}N'_\eps N_{\eps 1}N_\eps}}$
and apply the superquadratic property \eqref{superquadratic}
to obtain the bound
\begin{equation}
  \nonumber
  \frac{\sigma'}{4 s^*} \, 
  \frac{\delta_\eps^{\,2}}{\sqrt{\eps}} \,
  \frac{|g'_{\eps 1} g'_\eps| |q_\eps|}
       {\sqrt{N'_{\eps 1}N'_\eps N_{\eps 1}N_\eps}} 
  \leq \frac{1}{\varrho} \, \frac{\delta_\eps^{\,4}}{\eps} 
       \frac{{g'}_{\eps 1}^{\,2} {g'}_\eps^{\,2}}
            {N'_{\eps 1} N'_\eps N_{\eps 1} N_\eps} \,
       r^*\!\left( \frac{\sigma}{4 s^*} \right)
       G_{\eps 1} G_\eps
       + \frac{\varrho}{\eps \delta_\eps^{\,2}} \,
         r\!\left( \frac{\sqrt{\eps} \delta_\eps q_\eps}
                        {G_{\eps 1} G_\eps} \right) 
         G_{\eps 1} G_\eps \,.
\end{equation}
We then argue as we did to prove assertion \eqref{term1} from
\eqref{term1-Bnd}.
\end{proof}



\section{Limit of the Boundary Mass-Flux Term}

   In this section we show that as $\eps\to0$ the boundary term 
vanishes in the weak form of the approximate local conservation of 
mass that is obtained by setting $\zeta=1$ in \eqref{weak-moment}.
This is the key step in establishing the limiting mass conservation
equation \eqref{lim-mass-bdry} from \eqref{weak-moment}, as the limit
for all the other terms are obtained exactly as they were when we
established the interior mass conservation equation \eqref{lim-mass}.
More specifically, we prove the following.

\begin{Prop}
\label{Boundary-Prop}
For every $[t_1,t_2]\subset[0,\infty)$ and every
$\chi\in C^1(\overline{\Omega})$ one has
\begin{equation}
  \nonumber
  \lim_{\eps\to0} 
  \int_{t_1}^{t_2} \!\! \int_{\del\Omega} 
  \chi \, \< v \, \gamma \gTld_\eps \> \DOT \nRm \, \dsig \, \dt 
  = 0 \,.
\end{equation}
\end{Prop}

\begin{proof}
Denote the boundary mass-flux term as
\begin{equation}
  \nonumber
  j_\eps = \int_{t_1}^{t_2} \!\! \int_{\del\Omega} 
           \chi \, \< v \, \gamma \gTld_\eps \> \DOT \nRm \, 
           \dsig \, \dt \,.
\end{equation}
The renormalized boundary condition \eqref{renorm-bc} can be 
expressed as
\begin{equation}
  \nonumber
  \gamma_- \gTld_\eps 
  = L \left( \frac{\gHat_\eps}
                  {1 + \delta_\eps^{\,2} \gHat_\eps^{\,2}} \right) \,,
  \quad \hbox{where} \quad 
  \gHat_\eps 
  = \gamma_+ \big( (1 - \alpha) g_\eps
                   + \alpha \< \gamma_+ g_\eps \>_{\pO} \big) \,. 
\end{equation}
It follows that 
\begin{equation}
  \nonumber
\begin{aligned}
  j_\eps
  & = \int_{t_1}^{t_2} \!\! \int_{\del\Omega} \chi \, 
      \big\< v \, (\gamma_+ \gTld_\eps
                   + \gamma_- \gTld_\eps) \big\>   
      \DOT \nRm \, \dsig \, \dt 
\\
  & = \int_{t_1}^{t_2} \!\! \int_{\del\Omega} \chi \, 
      \left\< v \bigg( \gamma_+ \gTld_\eps
                       - \frac{\gHat_\eps}
                              {1 + \delta_\eps^{\,2} \gHat_\eps^{\,2}}
                       \bigg) \right\> 
      \DOT \nRm \, \dsig \, \dt 
\\
  & = \alpha \int_{t_1}^{t_2} \!\! \int_{\del\Omega} \chi \, 
      \left\< v \, \gamma_+ 
              \frac{(g_\eps - \< \gamma_+ g_\eps \>_{\pO})
                    (1 - \delta_\eps^{\,2} g_\eps \gHat_\eps)}
                   {(1 + \delta_\eps^{\,2} g_\eps^{\,2})
                    (1 + \delta_\eps^{\,2} \gHat_\eps^{\,2})} 
              \right\> \DOT  \nRm \, \dsig \, \dt \,.
\end{aligned}
\end{equation}
If $\alpha=0$ then we are done.  If $\alpha>0$ then set
$\gamma_\eps=\gamma_+\big(g_\eps-\<\gamma_+g_\eps\>_{\pO}\big)$, so 
that
\begin{equation}
  \label{j-gamma}
  j_\eps
  = \alpha \int_{t_1}^{t_2} \!\! \int_{\del\Omega} \chi 
    \int_{v\DOT\nRm>0} 
        \gamma_\eps \,  
        \frac{1 - \delta_\eps^{\,2} \gamp g_\eps \, \gHat_\eps}
             {(1 + \delta_\eps^{\,2} \gamp g_\eps^{\,2})
              (1 + \delta_\eps^{\,2} \gHat_\eps^{\,2})} \,
        v \DOT n \, M \dv \, \dsig \, \dt \,.
\end{equation}
The idea will now be to control $\gamma_\eps$ with the bound on the 
Darroz\`{e}s-Guiraud information \eqref{DG-inform} given by the
entropy inequality \eqref{entropy-inequality} and entropy bound
\eqref{entropy-bound}.   More specifically, we will show in Lemma 
\ref{gamma-Lem} that 
\begin{equation}
  \label{gamma-lim-1}
  \lim_{\eps\to0}
  \int_{t_1}^{t_2} \!\! \int_{\del\Omega} \chi \int_{v\DOT\nRm>0}
      \frac{\gamma_\eps}
           {(1 + \delta_\eps^{\,2} \gamp g_\eps^{\,2})
            (1 + \delta_\eps^{\,2} \gHat_\eps^2)} \,
      v \DOT \nRm \, M \dv \, \dsig \, \dt = 0 \,,
\end{equation}
and that
\begin{equation}
  \label{gamma-lim-2}
  \lim_{\eps\to0}
  \int_{t_1}^{t_2} \!\! \int_{\del\Omega} \chi \int_{v\DOT\nRm>0}
      \gamma_\eps \,  
      \frac{\delta_\eps^{\,2} \gamp g_\eps \, \gHat_\eps}
           {(1 + \delta_\eps^{\,2} \gamp g_\eps^{\,2})
           (1 + \delta_\eps^{\,2} \gHat_\eps^{\,2})} \,
      v \DOT n \, M \dv \, \dsig \, \dt = 0 \,.
\end{equation}
Propostion \ref{Boundary-Prop} for the case $\alpha>0$ will 
then follow from (\ref{j-gamma}-\ref{gamma-lim-2}) upon 
proving Lemma \ref{gamma-Lem}.
\end{proof}

\smallskip

\begin{Lem}
\label{gamma-Lem}
Let $\alpha>0$.  Then the limits \eqref{gamma-lim-1} and
\eqref{gamma-lim-2} hold.
\end{Lem}   

\begin{proof}
Following \cite{MS}, we employ the decomposition 
\begin{equation}
  \label{gamma-Decomp}
  \gamma_\eps = \gamma_\eps^{(1)} + \gamma_\eps^{(2)} \,, \quad
  \hbox{where} \quad
  \gamma_\eps^{(1)}
  = \gamma_\eps
    \One_{\{\gamp G_\eps 
            \leq 2 \langle \gamp G_\eps\rangle_{\pO}
            \leq 4 \gamp G_\eps \}} \,.
\end{equation}
By arguing as in Lemma 6.1 of \cite{MS} extended to our more general
class of collision kernels as in Lemma 6 of \cite{JM}, we obtain 
\begin{align}
  \label{bound-1}
  \frac{\gamma_\eps^{(1)}}
       {(1 + \delta_\eps^{\,2} \gamp g_\eps^2)^{1/4}}
  \quad & \mbox{is bounded in} \quad
  L^2_{loc}(\dt;L^2(|v\DOT\nRm|M\dv\,\dsig)) \,,
\\
  \label{bound-2}
  \frac{\gamma_\eps^{(1)}}
       {(1 + \delta_\eps^{\,2} \<\gamp g_\eps\>_{\pO}^{\,2})^{1/4}}
  \quad & \mbox{is bounded in} \quad
  L^2_{loc}(\dt;L^2(|v\DOT\nRm|M\dv\,\dsig)) \,,
\\
  \label{bound-3}
  \frac{1}{\delta_\eps} \gamma_\eps^{(2)} 
  \quad & \mbox{is bounded in} \quad 
  L^1_{loc}(\dt;L^1(|v\DOT\nRm|M\dv\,\dsig)) \,.
\end{align}

   To prove \eqref{gamma-lim-1} we use the fact that
$\<\gamma_\eps\>_{\pO}=0$ and the decompositon \eqref{gamma-Decomp} 
to write 
\begin{equation}
  \nonumber
\begin{aligned} 
  & \int_{t_1}^{t_2} \!\! \int_{\del\Omega} \chi \int_{v\DOT\nRm>0}
        \frac{\gamma_\eps}
             {(1 + \delta_\eps^{\,2} \gamp g_\eps^{\,2})
              (1 + \delta_\eps^{\,2} \gHat_\eps^{\,2})} \,
        v \DOT \nRm \, M \dv \, \dsig \, \dt
\\
  & = \int_{t_1}^{t_2} \!\! \int_{\del\Omega} \chi \int_{v\DOT\nRm>0}
      \bigg( \frac{\gamma_\eps}
             {(1 + \delta_\eps^{\,2} \gamp g_\eps^{\,2})
              (1 + \delta_\eps^{\,2} \gHat_\eps^{\,2})}
             - \frac{\gamma_\eps}
                    {{(1 + \delta_\eps^{\,2} \< \gamp g_\eps \>_{\pO}^{\,2})^2}}
             \bigg) v \DOT \nRm \, M \dv \, \dsig \, \dt
\\
  & = \int_{t_1}^{t_2} \!\! \int_{\del\Omega} \chi \int_{v\DOT\nRm>0}
      \bigg( \frac{\gamma_\eps^{(1)}}
             {(1 + \delta_\eps^{\,2} \gamp g_\eps^{\,2})
              (1 + \delta_\eps^{\,2} \gHat_\eps^{\,2})}
             - \frac{\gamma_\eps^{(1)}}
                    {{(1 + \delta_\eps^{\,2} \< \gamp g_\eps \>_{\pO}^{\,2})^2}}
             \bigg) v \DOT \nRm \, M \dv \, \dsig \, \dt
\\
  & \quad \,
      + \int_{t_1}^{t_2} \!\! \int_{\del\Omega} \chi \int_{v\DOT\nRm>0}
            \frac{\gamma_\eps^{(2)}}
                 {(1 + \delta_\eps^{\,2} \gamp g_\eps^{\,2})
                  (1 + \delta_\eps^{\,2} \gHat_\eps^{\,2})} \,
            v \DOT \nRm \, M \dv \, \dsig \, \dt
\\
  & \quad \,
      - \int_{t_1}^{t_2} \!\! \int_{\del\Omega} \chi \int_{v\DOT\nRm>0}
            \frac{\gamma_\eps^{(2)}}
                 {{(1 + \delta_\eps^{\,2} \< \gamp g_\eps \>_{\pO}^{\,2})^2}} \,
            v \DOT \nRm \, M \dv \, \dsig \, \dt \,.
\end{aligned}
\end{equation}
The last two terms on the right-hand side above vanish as $\eps\to0$
by the bound \eqref{bound-3}.  

   To show that the first term on the right-hand side above also 
vanishes as $\eps\to0$, we observe from the bounds \eqref{bound-1} 
and \eqref{bound-2}, the two terms in its integrand are relatively 
compact in $\wL^2_{loc}(\dt;\wL^2(|v\DOT\nRm|M\dv\,\dsig))$.  
Their difference is
\begin{equation}
  \nonumber
\begin{aligned}
  & \frac{\gamma^{(1)}_\eps}
         {(1 + \delta_\eps^{\,2} \gamp g_\eps^{\,2})
          (1 + \delta_\eps^{\,2} \gHat_\eps^{\,2})}
    - \frac{\gamma^{(1)}_\eps}
           {(1 + \delta_\eps^{\,2} \< \gamp g_\eps \>_{\pO}^{\,2})^2}
\\
  & = \gamma^{(1)}_\eps
      \frac{[(1 + \delta_\eps^{\,2} \< \gamp g_\eps \>_{\pO}^{\,2})^2
             - (1 + \delta_\eps^{\,2} \gamp g_\eps^{\,2})^2]
            + [(1 + \delta_\eps^{\,2} \gamp g_\eps^{\,2}) 
               - (1 + \delta_\eps^{\,2} \gamp g_\eps^{\,2}) 
                 (1 + \delta_\eps^{\,2} \gHat_\eps^{\,2})]}
           {(1 + \delta_\eps^{\,2} \gamp g_\eps^{\,2})
            (1 + \delta_\eps^{\,2} \gHat_\eps^{\,2})
            (1 + \delta_\eps^{\,2} \< \gamp g_\eps \>_{\pO}^2)^2} \,.
\end{aligned}
\end{equation}
Noting that $\gamp g_\eps - \gHat_\eps = \alpha \gamma_\eps$, 
we see that
\begin{equation}
  \label{gamma1-difference1}
\begin{aligned}
  & \frac{\gamma^{(1)}_\eps}
         {(1 + \delta_\eps^{\,2} \gamp g_\eps^{\,2})
          (1 + \delta_\eps^{\,2} \gHat_\eps^{\,2})}
    - \frac{\gamma^{(1)}_\eps}
           {(1 + \delta_\eps^{\,2} \< \gamp g_\eps \>_{\pO}^{\,2})^2}
\\
  & = \delta_\eps
      \left( \frac{\gamma^{(1)}}
                  {(1 + \delta_\eps^{\,2} \gamma g_\eps^{\,2})^{1/4}} \right)^2 
      \One_{\{\gamp G_\eps 
              \leq 2 \langle \gamp G_\eps \rangle_{\pO}
              \leq 4 \gamp G_\eps\}} 
      (S^{(1)}_\eps + S^{(2)}_\eps) \,,
\end{aligned}
\end{equation}
where
\begin{equation}
  \nonumber
\begin{aligned}
  S^{(1)}_\eps
  & = \frac{\delta_\eps (\gamp g_\eps + \< \gamp g_\eps \>_{\pO})
            (2 + \delta_\eps^{\,2} \gamp g_\eps^{\,2} 
               + \delta_\eps^{\,2} \< \gamp g_\eps \>_{\pO}^{\,2})}
           {(1 + \delta_\eps^{\,2} \gamp g_\eps^{\,2})^{1/2} 
            (1 + \delta_\eps^{\,2} \gHat_\eps^{\,2})
            (1 + \delta_\eps^{\,2} \< \gamp g_\eps \>_{\pO}^{\,2})^2} \,,
\\
  S^{(2)}_\eps
  & = - \frac{\alpha \delta_\eps (\gamp g_\eps + \gHat_\eps)
              (1 + \delta_\eps^{\,2} \gamp g_\eps^{\,2})}
             {(1 + \delta_\eps^{\,2} \gamp g_\eps^{\,2})^{1/2}
              (1 + \delta_\eps^{\,2} \gHat_\eps^{\,2})
              (1 + \delta_\eps^{\,2} \<\gamp g_\eps\>_{\pO}^{\,2})^2} \,.
\end{aligned}
\end{equation}
Note that $\gamp G_\eps\leq 2\< \gamp G_\eps \>_{\pO}\leq 4\gamp G_\eps$ 
implies that
\begin{equation}
  \nonumber
  \delta_\eps \gamp g_\eps 
  \leq 2 \delta_\eps \< \gamp g_\eps \>_{\pO} + 1 \,,
  \quad \mbox{and} \quad
  \delta_\eps \< \gamp g_\eps \>_{\pO} 
  \leq 2 \delta_\eps \gamp g_\eps + 1 \,.
\end{equation}
It is easy to see that both
\begin{equation}
  \nonumber
  \mathbf{1}_{\{\gamp G_\eps \leq 2 \< \gamp G_\eps \>_{\pO}
                \leq 4 \gamp G_\eps\}} S^{(1)}_\eps
  \quad \mbox{and} \quad
  \mathbf{1}_{\{\gamp G_\eps \leq 2 \< \gamp G_\eps \>_{\pO}
                \leq 4 \gamp G_\eps\}} S^{(2)}_\eps
\end{equation}
are bounded in $L^\infty$.  Noting that the extra $\delta_\eps$ in 
front of \eqref{gamma1-difference1}, we conclude that for any 
convergent subsequence,
\begin{equation}
  \nonumber
  \frac{\gamma^{(1)}_\eps}
       {(1 + \delta_\eps^{\,2} \gamp g_\eps^{\,2})
        (1 + \delta_\eps^{\,2} \gHat_\eps^{\,2})}
        -\frac{\gamma^{(1)}_\eps}
              {(1 + \delta_\eps^{\,2} \< \gamp g_\eps \>_{\pO}^{\,2})^2} 
  \to 0 \,,
  \quad \mbox{in} \quad 
  \wL^2_{loc}(\dt;\wL^2(|v\DOT\nRm|M\dv\,\dsig)) \,,
\end{equation}
as $\eps\to0$.  This establishes limit \eqref{gamma-lim-1}.

   To prove limit \eqref{gamma-lim-2} separate 
$\gamma_\eps = \gamma^{(1)}_\eps + \gamma^{(2)}_\eps$, using the 
bound \eqref{bound-3} for $\gamma^{(2)}_\eps$, it is easy to estimate 
that
\begin{equation}
  \nonumber
  \delta_\eps
  \left\| \frac{\alpha}{\delta_\eps} \gamma^{(2)}_\eps 
          \right\|_{L^1_{loc}(\dt;L^1(|v\DOT\nRm|M\dv\,\dsig))}
  \frac{|\delta_\eps \gamp g_\eps \, \delta_\eps \gamp \gHat_\eps|}
       {(1 + \delta_\eps^{\,2} \gamp g_\eps^{\,2}) 
        (1 + \delta_\eps^{\,2} \gHat_\eps^{\,2})}
  \leq C \tfrac{1}{4} \delta_\eps \,.
\end{equation}
For the $\gamma^{(1)}_\eps$ part, from the $L^2$ bounds 
\eqref{bound-1} and \eqref{bound-2},
\begin{equation}
  \label{product-1}
  \frac{\gamma^{(1)}_\eps}
       {\sqrt{1 + \delta_\eps^{\,2} \gamp g_\eps^{\,2}}}
\end{equation}
is relatively compact in
$\wL^1_{loc}(\dt;\wL^1(|v\DOT\nn|M\dv\,\dsig))$.  Use
the fact that
\begin{equation}
  \label{product-2}
  \frac{\delta_\eps \gamp g_\eps \, \delta_\eps \gHat_\eps}
       {\sqrt{1 + \delta_\eps^{\,2} \gamp g_\eps^{\,2}} 
        (1 + \delta_\eps^{\,2} \gHat_\eps^{\,2})}
\end{equation}
is bounded in $L^\infty$ and goes to 0 a.e.  Then again by the
Product Limit Theorem of \cite{BGL3}, the product of
\eqref{product-1} and \eqref{product-2} goes to 0 in
$L^1_{loc}(\dt;L^1(|v\DOT\nRm|M\dv\,\dsig))$ as $\eps\to0$.  
We thereby finish the proof of the Lemma.
\end{proof}

\smallskip
\noindent
{\bf Remark.} The most important difference between the acoustic limit 
and the incompressible limits (Stokes in \cite{MS} and Navier-Stokes 
in \cite{JM}) is that the compactness of the renormalized traces 
$\gamma \gTld_\eps$ in the acoustic limit case is not available. 
The pointwise convergence $\delta_\eps \gTld_\eps \to 0$ a.e. is 
also unavailable. (compare Lemma 5.2 in \cite{MS}.)  In contrast, for the 
incompressible limits the entropy bounds from boundary provide 
{\em a priori} estimates on the quantity 
$\gamma_\eps = \gamp g_\eps - \mathbf{1}_{\Sigma_+\<\gamp g_\eps\>_{\pO}}$. 
Specifically, we have the $L^2$ bound on 
$\frac{1}{\delta_\eps}\frac{\gamma^{(1)}_\eps}{n_\eps}$ with some 
renormalizer $n_\eps$, see bounds (6.2) and (6.3) in \cite{MS}. 
However, in the acoustic limit, because of the acoustic scaling, we have 
only the $L^2$ bound on $\frac{\gamma^{(1)}_\eps}{n_\eps}$ which is much 
weaker than in the incompressible limits cases.

\end{document}